\documentclass[a4paper]{article}

\usepackage{amsmath,amssymb,amsthm}
\usepackage[graph]{xy}
\theoremstyle{plain}      
    \newtheorem{theorem}{Theorem}[section]

\theoremstyle{definition}
    
    \newtheorem{example}[theorem]{Example}
\theoremstyle{remark}
    \newtheorem*{remark}{Remark}

\DeclareMathOperator{\Wien}{Wien}

\newcommand{\A}{\ensuremath{\mathcal{A}}}
\newcommand{\B}{\ensuremath{\mathcal{B}}}
\newcommand{\C}{\ensuremath{\mathcal{C}}}
\renewcommand{\P}{\ensuremath{\mathcal{P}}}
\newcommand{\X}{\ensuremath{\mathcal{X}}}
\newcommand{\V}{\ensuremath{\mathcal{V}}}
\newcommand{\Z}{\ensuremath{\mathcal{Z}}}
\newcommand{\Vect}{\ensuremath{\mathbf{Vect}}}
\newcommand{\Set}{\ensuremath{\mathbf{Set}}}
\newcommand{\Cat}{\ensuremath{\mathbf{Cat}}}
\newcommand{\fd}{\ensuremath{\mathrm{fd}}}
\newcommand{\op}{\ensuremath{\mathrm{op}}}
\newcommand{\ox}{\ensuremath{\otimes}}

\newcommand{\x}{\ensuremath{\times}}
\newcommand{\ra}{\ensuremath{\rightarrow}}
\newcommand{\xra}{\ensuremath{\xrightarrow}}

\begin{document}

\title{Association schemes, classical RCFT's, and centres of monoidal functor
categories}
\author{Brian Day}

\maketitle

\begin{abstract}
Here we describe three straightforward examples of what was called a graphic
Fourier transformation in~\cite{4}. At least two of these examples may be
viewed simply as monoidal comonads on suitable monoidal closed functor
categories, but the third example, which involves ``centres'' of monoidal
closed functor categories, is generally not comonadic. For the first two
examples (i.e., association schemes and RCFT's), a more elaborate
``probicategory'' set-up was envisaged in an earlier version of this note, but
many readers missed the main point so it is simplified (hopefully) below.
\end{abstract}

\section{Introduction}

A functor category such as $[\A,\Vect]$ or $[\A,\Set]$ behaves naively like
a function algebra, both types of structure being familiar manifestations of
``form'' and ``function''.

Upon setting $\V = \Vect$ (or $\Set$) and taking $\A$ and $\X$ to be two
(small) promonoidal $\V$-categories, we let $[\A,\V]$ and $[\X,\V]$ denote
the corresponding convolution functor categories into $\V$ (as in~\cite{2}).
A $\V$-functor
\[
    \hat{K} : [\A,\V] \ra [\X,\V]
\]
is then called a ``graphic'' Fourier transformation~\cite{4} if

\begin{enumerate}
\item $\hat{K}$ is multiplicative (i.e., $\hat{K}$ preserves the tensor
product and tensor unit of $[\A,\V]$ up to natural isomorphism),

\item $\hat{K}$ is conservative (i.e., $\hat{K}$ reflects isomorphisms),

\item $\hat{K}$ is cocontinuous (i.e., $\hat{K}$ has a right $\V$-adjoint,
denoted here by $\check{K}$).
\end{enumerate}
Sometimes the convolutions $[\A,\V]$ and $[\X,\V]$ both have (contravariant)
involutions on them and we ask in addition that $\hat{K}$ should preserve
this.

Below we mention two elementary examples in which (in the first instance) the
category $\A$ is discrete. We also describe a third type of example based on
the ``centres'' of certain monoidal closed functor categories into $\V$.

For the general theory, the ``image'' of $\hat{K}$ -- here called the Wiener
category (cf.~\cite{4} \S 1.3) and denoted by $\Wien(K)$ -- was discussed
briefly in~\cite{4}. Roughly speaking, its objects are the functors of the form
$\hat{K}(f)$ for some $f \in [\A,\V]$, and its maps are precisely those
natural transformations
\[
    \alpha : \hat{K}(f) \ra \hat{K}(g)
\]
in $[\X,\V]$ for which
\[
    \hat{K}\check{K}(\alpha) \hat{K}(\eta) = \hat{K}(\eta)\alpha,
\]
where $\eta$ denotes the (monomorphic) unit of the basic adjunction $\hat{K}
\dashv \check{K}$. These so-called ``regular'' morphisms $\alpha \in [\X,\V]$
have been further characterized in many special examples in the literature.

Moreover, the property that $\hat{K}$ should be \emph{multiplicative} amounts,
under the left Kan extension process, to the existence of two natural
isomorphisms:
\begin{align*}
    \int^{y,z \in \X} K(a,y) \ox K(b,z) \ox p(y,z,x) & \xra{\cong} 
    \int^{c \in \A} K(c,x) \ox p(a,b,c) \\
    j(x)  & \xra{\cong} \int^{c \in \A} K(c,x) \ox j(c)
\end{align*}
in $\V$, where the functor
\[
    K : \A^\op \ox \X \ra \V
\]
denotes the tensor-hom transpose in $\V$-$\Cat$ of the composite
\[
    \A^\op \xra{\text{Yoneda}} [\A,\V] \xra{\hat{K}} [\X,\V],
\]
and was called the (multiplicative) \emph{kernel} of $\hat{K}$ in~\cite{4}.

Then we obtain the monoidal equivalence of categories
\[
    \hat{K} : [\A,\V] \ra \Wien(K).
\]

\section{Terminology}

A $\Set$-promonoidal category $(\A,n,i)$ where $\A$ is a groupoid with antipode
$S:\A^\op \ra \A$ ($S^2 \cong 1$) is called \emph{($S$-)precompact} if both
\begin{align*}
    n(a,b,c) &\cong n(Sb, Sa, Sc) \text{ and} \\
    i(c) &\cong i(Sc)
\end{align*}
dinaturally in $a,b,c \in \A$. Then, if $n$ and $i$ are finite, we have both
\begin{align*}
    k[n](a,b,c)^* &\cong k[n](Sb, Sa, Sc) \text{ and} \\
    k[i](c)^* &\cong k[i](Sc)
\end{align*}
naturally in $a,b,c \in \A$, where $k$ is any field. Thus, we call any finite
$\Vect_k$-promonoidal category $(A,p,j)$, with a $k$-linear antipode $S$ ($S^2
\cong 1$), \emph{precompact} if both
\begin{align*}
    p(a,b,c)^* &\cong p(Sb, Sa, Sc) \text{ and} \\
    j(c)^* &\cong j(Sc)
\end{align*}
naturally in $a,b,c \in \A$.

Note that, when the above data is available we usually have the ``structure''
maps
\begin{align*}
    f^* \ox g^* &\ra (g \ox f)^* \\
    I &\ra I^*,
\end{align*}
where $f^*$ is defined by $f^*(a) = f(Sa)^*$, giving $(-)^*$ the structure of a
contravariant monoidal functor on $[\A,\Vect_\fd]$. This is so, for example,
when $\A$ is ``Frobenius'' (see~\cite{5} --- there is a family of maps
\[
    \delta : \A(a,b) \ra \A(c,b) \ox \A(a,c)
\]
in $\V$ which is natural in $a,b \in \A$). In this case there is a ``comparison
of integrals''
\[
    \int^a h(a,a) \ra \int_b h(b,b)
\]
for each $k$-linear functor $h : \A^\op \ox \A \ra \Vect_k$ which is derived
directly from the composite natural transformation
\[
    h(a,a) \xra{\text{``}\delta\text{''}} \A(a,b) \ox h(b,a)
    \xra{\text{``}\mu\text{''}} h(b,b).
\]
Thus one has the composite:
\begin{align*}
    (f^* \ox g^*)(c) &:= \int^{a,b} f(Sa)^* \ox g(Sb)^* \ox p(a,b,c)
        & \text{(by definition of $f^* \ox g^*$)} \\
    &\cong \int^{a,b} f(a)^* \ox g(b)^* \ox p(Sa,Sb,c) & (S^2 \cong 1) \\
    &\ra \int_{a,b} g(b)^* \ox f(a)^* \ox p(Sa,Sb,c)
        & \text{(by ``comparison of integrals'')} \\
    &\cong \left( \int^{a,b} g(b) \ox f(a) \ox p(b,a,Sc)\right)^*
        & \text{(when $\A$ is precompact)} \\
    &= (g \ox f)(Sc)^* & \text{(by definition of $g \ox f$)} \\
    &= (g \ox f)^*(c).
\end{align*}

\section{Discrete association schemes}

Let $X$ be a set and let $S \subset \P(X \x X)$ be a fixed partition of $X \x
X$ such that $\Delta \in S$ and $s^* \in S$ when $s \in S$, where $\Delta
\subset X \x X$ is the diagonal set, and $s^*$ is the reverse relation to $s$.
The relevant kernel
\[
    K : S \x X \x X \ra \Set
\]
is given here by
\[
    K(s,x,y) = s(x,y) =
        \begin{cases}
            1 & \text{if $(x,y) \in s$} \\ 
            0 & \text{otherwise}.
        \end{cases}
\]
The resulting transformation
\[
    \hat{K} : [S,\Set] \ra [X \x X, \Set]
\]
then maps $f$ to
\begin{align*}
    \hat{K}(f)(x,y) &= \sum_s K(s,x,y) \x f(s) \\
                    &= \sum_{(x,y) \in s} f(s).
\end{align*}

A (discrete) \emph{association scheme} on $X$ (after~\cite{10}, for example)
consists of the given partition $S$ and a map
\[
    N : S \x S \x S \ra \Set
\]
together with a bijection
\[
    \sum_r N(s,t,r) \x r \xra{\cong} s \circ t
\]
for each $(s,t) \in S \x S$, where $s \circ t$ denotes the $\Set$-matrix
composite
\[
    s \circ t(x,y) = \sum_{z \in X} s(x,z) \x t(z,y).
\]
It follows readily from the associativity of this composition that there is a
proassociativity bijection
\[
    \alpha : \sum_u N(s,t,u) \x N(u,r,v) \xra{\cong} \sum_u N(s,u,v) \x
    N(t,r,u)
\]
for all $s,t,r,v \in S$ since $\sum_v M(v) \x v \cong \sum_v P(v) \x v$ implies
$M(v) \cong P(v)$ for all $v \in S$ by the definition of ``partition''.
Together with the diagonal ``object'' $\Delta$, this data represents a
discrete promonoidal structure called $(S,N,J)$ on the ``category'' $S$ (where
the prounit $J:S \ra \Set$ is given by
\[
    J(s) = \begin{cases} 1 & \text{if $s = \Delta$} \\ 0 & \text{otherwise}.
           \end{cases}
\]
Moreover, we may suppose (see~\cite{10}) that 
\[
    N(s,t,r) = N(t^*,s^*,r^*)
\]
so that $(S,N,J)$ is what we have called precompact in \S2. Additionally, if
also $X$ is finite and
\[
    N(s,t,r^*) = N(t,r,s^*),
\]
then the monoidal closed convolution structure on $[S,\Vect_\fd]$ is in fact
\emph{compact}, as shown in the arXiv note~\cite{3}.

Now it is easily seen that the transformation
\[
    \hat{K} : [S,\Set] \ra [X \x X, \Set]
\]
is multiplicative (where $[S,\Set]$ has the convolution structure~\cite{2} from
$(S,N,J)$), conservative, and cocontinuous. Also, if we set $f^*(s) = f(s^*)$
then we get
\begin{align*}
    \hat{K}(f^*)(x,y) &= \sum_{(x,y) \in s} f^*(s) \\
                      &= \sum_{(x,y) \in s} f(s^*) \\
                      &= \sum_{(y,x) \in s^*} f(s^*) \\
                      &= \sum_{(y,x) \in t} f(t) \\
                      &= \hat{K}(f)^*(x,y).
\end{align*}
Thus, $\hat{K}$ becomes a graphic Fourier transformation on $[S,\Set]$, in the
sense of~\cite{4}.

In the case where the set $X$ is \emph{finite}, we can further view each
element $s \in S$ as a ``$k$-bimodule'' $M_s$, for the given field $k$, and
write
\[
    M_s \circ M_t \cong \bigoplus_r N(s,t,r) \cdot M_r
\]
where ``$\circ$'' denotes matrix composition, and ``$\cdot$'' denotes copowers
in $[X \x X, \Vect_\fd]$. The functor
\[
    \hat{K} : [S, \Vect_\fd] \ra [X \x X, \Vect_\fd]
\]
given by
\[
    \hat{K}(f)(x,y) = \bigoplus_{(x,y) \in s} f(s)
\]
then becomes a $k$-linear graphic Fourier transformation on $[S, \Vect_\fd]$.
Here we can actually set
\[
    f^*(s) := f(s^*)^*
\]
and obtain
\[
    \hat{K}(f^*)(x,y) \cong \hat{K}(f)^*(x,y)
\]
for all $(x,y) \in X \x X$; i.e.,
\[
    \hat{K}(f^*) \cong \hat{K}(f)^*.
\]

\section{Rational conformal field theories}

A discrete RCFT consists of a finite (discrete) set $S$, containing a
distinguished element $0$, and a set of finite sets $N(x,y,z)$ indexed by
elements $x,y,z$ of $S$. We suppose also that there are proassociativity and
prounit maps (bijections)
\begin{align*}
    \alpha: \sum_u N(x,y,u) \x N(u,z,v) &\xra{\cong}
        \sum_u N(x,u,v) \x N(y,z,u) \\
        \lambda : N(0,y,z) &\xra{\cong} S(y,z) \\
        \rho : N(x,0,z) &\xra{\cong} S(x,z),
\end{align*}
which satisfy some standard coherence conditions. There is usually assumed to
be an involution $(-)^*$ on $S$ satisfying the cyclic relation
\[
    N(x,y,z^*) \cong N(y,z,x^*),
\]
together with a ``braiding''
\[
    N(x,y,z) \cong N(y,x,z).
\]
(See Moore and Sieberg~\cite{8} for details.)

Thus we have a $k$-linear functor
\[
    p : \A^\op \ox \A^\op \ox \A \ra \Vect_\fd,
\]
where $\A$ denotes the discrete $\Vect$-category on $S$, by considering the
free $k$-vector spaces on the relevant finite sets, that is
\[
    p(x,y,z) = k[N](x,y,z).
\]
Together with the corresponding induced proassociativity and prounit
isomorphisms, we obtain a braided (finite) $*$-autonomous promonoidal
$\Vect$-category $(\A,p,j)$. We could clearly generalize this situation
considerably (in theory at least). For example, the arXiv note~\cite{3} contains
further details of compact convolution categories $[\A,\Vect_\fd]$ in the case
when the domain category $\A$ is more than just discrete or finite.

Always the representation (or ``Cayley'') functor
\[
    \hat{K} : [\A,\Vect_\fd] \ra [\A^\op \ox \A,\Vect_\fd],
\]
which maps $f$ to
\[
    \hat{K}(f)(y,z) = \int^x f(x) \ox p(x,y,z),
\]
can be seen to be multiplicative, conservative, and cocontinuous (as is
generally the case), thus providing examples of $k$-linear graphic Fourier
transforms on certain compact convolutions of the form $[\A,\Vect_\fd]$.
Indeed, if we suppose further that $(\A,p,j)$ is any closed category equipped
with a suitable involution
\[
    (-)^* : \A^\op \ra \A
\]
(that is, $p(x,y,z) \cong \A(x,[y,z])$ with $[x,y]^* \cong [y,x]$), then we
also have
\begin{align*}
    \hat{K}(f^*)(x,y) &= \int^z f^*(z) \ox \A(z,[x,y]) \\
        &\cong f^*[x,y] \\
        &\cong f[y,x]^* \\
        &\cong \hat{K}(f)(y,x)^* \\
        &= \hat{K}(f)^*(x,y),
\end{align*}
thus completing the list of properties for a graphic transformation listed in
the introduction.

\section{Centres of some monoidal functor categories}

If $(\A,p,j)$ is any small promonoidal $\V$-category, where $\V$ is either
$\Vect$ or $\Set$, then the centre $\Z[\A,\V]$ of the monoidal convolution
$[\A,\V]$ has the property that the underlying-object functor
\[
    U : \Z[\A,\V] \ra [\A,\V]
\]
is multiplicative, conservative, and cocontinuous. Consequently, if we can show
that a particular centre $\Z[\A,\V]$ is monoidally equivalent to a complete
convolution functor category of the form $[\B,\V]$, then the composite functor
\[
    \xygraph{{[\B,\V]}="1" [r(2)]{[\A,\V]}="2" [l(1)d]{\Z[\A,\V]}="3"
        "1":"2" ^-{\hat{K}}
        "1":"3"
        "3":"2" _-U}
\]
becomes a graphic Fourier transformation.

We remark that some efforts were made in~\cite{6} and~\cite{7} to find suitable
representing categories $\B$ for $\Z[\A,\V]$ and, in particular, the functor
\[
    \tilde{\Psi} : \Z[\A,\V] \ra [\Z(\A),\V]
\]
defined in~\cite{6}\S 2 is an equivalence of categories if and only if both
\begin{enumerate}
\item the canonical full embedding denoted by
\[
    \Psi:\Z(\A)^\op \subset \Z[\A,\V]
\]
in~\cite{6} is dense, and

\item existential quantification
\[
    \exists_U : [\Z(\A),\V] \ra [\A,\V]
\]
is conservative.
\end{enumerate}
Here $U:\Z(\A) \ra \A$ denotes the forgetful functor from the ``promagmal''
centre $\Z(\A)$ (see~\cite{6}\S 2) of the given promonoidal structure
$(\A,p,j)$; thus $\Z(\A)$ reproduces the usual braided monoidal centre of $\A$
when the promonoidal structure on $\A$ is in fact monoidal.

In this regard, we point out that, in the case of $\V = \Vect$, the full
embedding
\[
    \Psi:\Z(\A)^\op \subset \Z[\A,\Vect]
\]
is dense whenever all epimorphisms split in $[\A,\Vect]$ and each of the
functors
\[
    p(a,b,-) : \A \ra \Vect, \quad a,b \in \A,
\]
is finitely presentable in $[\A,\Vect]$; for example, if $\A$ is monoidal then
each functor
\[
    p(a,b,-) \cong \A(a \ox b,-) : \A \ra \Vect
\]
is certainly finitely presentable in $[\A,\Vect]$

\begin{example}
(cf.~\cite{9}) All epimorphisms split in the functor category
$[\text{Span}(\C),\Vect]$ when $\C =$ finite $G$-sets for a finite group $G$,
so that
\[
    \Psi : \Z(k_* \text{Span}(\C))^\op \subset \Z[k_* \text{Span}(\C),\Vect]
\]
is dense in this case, because $k_* \text{Span}(\C)$ (which denotes the
$k$-linearization of $\text{Span}(\C)$) is monoidal. Moreover,
\[
    \tilde{\Psi} : \Z[k_* \text{Span}(\C),\Vect] \simeq
    [\Z(k_* \text{Span}(\C)),\Vect]
\]
is an equivalence, which is not as refined as the main Theorem of~\cite{9} p.
4019, but is much more easily obtained.
\end{example}

\begin{example}
(cf.~\cite{1}) All epimorphisms split in the functor category $[\A,\Vect]$ for
any compact monoidal fusion category $\A$ over a field $k$ of suitable
characteristic, so here also $\Psi$ is dense, and
\[
    \tilde{\Psi} : \Z[\A,\Vect] \simeq [\Z(\A),\Vect]
\]
is an equivalence of monoidal categories.
\end{example}

The duality and precompactness conditions mentioned in \S 2 can also be
considered in this context.

\begin{remark}
We have considered just the \emph{centres} of monoidal functor categories here.
In fact, the results used above (from~\cite{6} and~\cite{7}) were originally
established in terms of ``lax'' centres, but there is no significant difference
in dealing with the centres.
\end{remark}



\bigskip
{\small\noindent
Department of Mathematics \\
Macquarie University \\
NSW, 2109, Australia}

\end{document}